\documentclass{llncs}

\usepackage{ucs}
\usepackage[utf8x]{inputenc}

\usepackage{pstricks,pst-plot,pst-node}
\usepackage{gensymb,enumerate,amssymb}

\title{Computable counter-examples to the Brouwer fixed-point theorem}
\author{Petrus H. Potgieter}
\institute{Department of Decision Sciences, University of South Africa (Pretoria), PO Box 392,  Unisarand,  0003,  Republic of South Africa, %fax:+27-12-429-4898, 
\email{php@member.ams.org}, \email{potgiph@unisa.ac.za}, \email{www.potgieter.org}}

\newcommand{\field}[1]{\ensuremath{\mathbb{#1}}}

\newcommand{\CR}{\ensuremath{\field{R}_c}}

\newcommand{\NZ}{\ensuremath{\field{N}_0}}

\begin{document}

\maketitle

\begin{abstract}
This paper is an overview of results that show the Brouwer fixed-point theorem (BFPT) to be essentially non-constructive and non-computable. The main results, the counter-examples of Orevkov and Baigger, imply that there is no procedure for finding the fixed point in general by giving an example of a computable function which does not fix any computable point. Research in reverse mathematics has shown the BFPT to be equivalent to the weak König lemma in RCA$_0$ (the system of recursive comprehension) and this result is illustrated by relating the weak König lemma directly to the Baigger example.
\end{abstract}

\section{Introduction}

We consider the Brouwer fixed-point theorem (BFPT) in the following form, where the standard unit interval is denoted by $I=[0,1]$.
\begin{theorem}[Brouwer]
Any continuous function $f:I^2 \rightarrow I^2$ has a fixed point, i.e. there exists an $x\in I^2$ such that $f(x)=x$.
\end{theorem}
A computable real number is a number for which a Turing machine exists that, on input $n$, produces a rational approximation with error no more than $2^{-n}$. A computable point is a point all the coordinates of which are computable reals. The notation
\begin{quote}
\begin{tabular}{rl}
\NZ & for the non-negative natural numbers;\\
\CR & for the set of computable reals;\\
$I_c$ & for $I \cap \CR$; and\\
$\delta X$ &  for the boundary of a set $X$, being $\overline{X} \cap \overline{X^c}$
\end{tabular}
\end{quote}
is also used. The two examples discussed use distinct definitions of a computable function of real variables.
\begin{description}
\item[Russian school] In the Russian school of Markov and others, a computable function maps computable reals to computable reals by a single algorithm for the function that translates an algorithm approximating the argument to an algorithm approximating the value of the functions. It need not be possible to extend a function that is computable in the Russian school to a continuous function on all of the reals. These functions are often called \textit{Markov-computable}.
\item[Polish school] In the Polish school of Lacombe, Grzegorczyk, Pour-El and Richards, and others, a function is computable on a region if it maps every every computable sequence of reals to a computable sequence of reals and it has a computable uniform modulus of continuity on the region  \cite{PourElRichards}.
\end{description}

\section{Orevkov's example for the Russian school}

One can construct a Markov-computable function $f$ through a computable mapping of descriptions of computable points $x \in I_c^2$ to descriptions of $f(x) \in I_c^2$, such that
$$f(x) \neq x \quad\quad \forall x \in I_c^2.$$
That is, no computable point is a fixed point for $f$. Unfortunately the $f$ which is constructed in this way, cannot be extended to a continuous function on $I^2$. This is the construction of \cite{Orevkov63}, another instance of which can be found in \cite{WR1999a}.

\begin{figure}
\centering
\psset{unit=3cm}
\begin{pspicture}[shift=-0.5](-0.25,0)(1.25,1.25)
\psframe[dimen=inner,linewidth=0pt](0,0)(1,1)
\pspolygon[linecolor=gray,fillstyle=solid,fillcolor=gray](0.5,0.2)(0.5,0.7)(0.8,0.7)(0.8,0.2)
\pspolygon[fillstyle=crosshatch](0.3,0.7)(0.3,0.9)(1.0,0.9)(1.0,0.7)
\pspolygon[fillstyle=crosshatch](0.3,0.1)(0.3,0.2)(1.0,0.2)(1.0,0.1)
%(1.0,0.2)(0.8,0.2)(0.8,0.7)
\end{pspicture}
$\cdots$
\begin{pspicture}[shift=-0.5](-0.25,0)(1.25,1.25)
\psframe[dimen=inner,linewidth=0pt](0,0)(1,1)
\pspolygon[fillstyle=crosshatch](0.3,0.7)(0.3,0.9)(1.0,0.9)(1.0,0.7)
\pspolygon[fillstyle=crosshatch](0.3,0.1)(0.3,0.2)(1.0,0.2)(1.0,0.1)
\pspolygon[linewidth=3pt,linecolor=gray](0.5,0.7)(0.5,0.7)(1.0,0.7)(1.0,0.7)
\pspolygon[linewidth=3pt,linecolor=gray](1.0,0.2)(1.0,0.7)(1.0,0.7)(1.0,0.2)
\pspolygon[linewidth=3pt,linecolor=gray](0.5,0.2)(0.5,0.2)(1.0,0.2)(1.0,0.2)
\end{pspicture}
\caption{Basic contraction in the Orevkov counter-example}
\label{fig:Ore}
\end{figure}
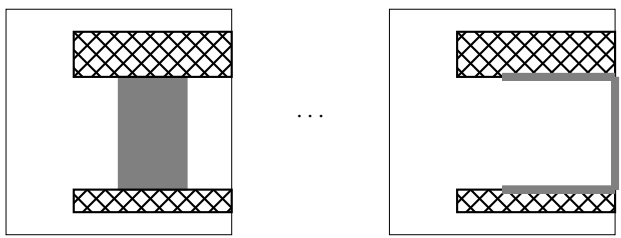

\begin{lemma}
Suppose $A_k$ is a sequence of rectangles in $I^2$ with computable vertices, disjoint interiors, and such that
\begin{enumerate}[(i)]
\item $\emptyset \neq \delta A_j \setminus \bigcup_{i < j} A_i$ for all $j$;
\item for each $j$ there exists $n > j$ such that $\delta A_j \subset \left( \bigcup_{i \leq n} A_i \right)\degree$; and
\item $I_c^2 \subseteq \bigcup_{i \geq 1} A_i$
\end{enumerate}
then there exists a Markov computable $g$, mapping $I_c^2$ to $\delta I_c^2$ and fixing $\delta I_c^2$.
\end{lemma}
The conditions ensure that 
\begin{enumerate}[(i)]
\item rectangles $A_j$, when added, have some part of their boundary in $I^2 \setminus \bigcup_{i < j} A_i$;
\item each $A_j$ is eventually closed off by new rectangles on all sides;
\item all computable points lie in $\bigcup_{i \geq 1} A_i$.
\end{enumerate}%
The function $f$ is obtained by composing $g$ with a $90\degree$ rotation. It therefore remains only to prove the lemma and the existence of a sequence of rectangles which is as required. Suppose that $g$ has been defined on $\bigcup_{i < j} A_i$. For
\begin{description}
\item[$\emptyset = \delta A_j \cap \bigcup_{i < j} A_i$:] let $g$ on $A_j$ consist of the simplest possible mapping to $\delta I^2$,
that fixes $\delta I^2$;
\item[$\emptyset \neq \delta A_j \cap \bigcup_{i < j} A_i$:] we extend $g$ to $A_j$  by using (i)---if $g$ has already been defined on the crosshatched set in Figure \ref{fig:Ore} then the definition can be extended to the solid gray set $A_j$ by composing a contraction of the solid gray set in Figure \ref{fig:Ore} to
$$\left( \delta A_j \cap \bigcup_{i < j} A_i\right) \cup \left( \delta A_j \cap \delta I^2 \right)$$
with the function $g$ as it has already been defined on the crosshatched set. This is always possible because, by construction of the $A_i$, our $A_j$ will always have at least \textit{two} sides non-contiguous with $\bigcup_{i < j} A_i$, at least one of which will not coincide with $\delta I^2$.
\end{description}
So far only condition (i) has been used. Conditions (ii) and (iii) are necessary for showing that $g$ is Markov-computable on $I^2$. Let a description of any $x \in I_c^2$ be given. We can find a description of $g(x)$ in the following way.
\begin{itemize}
\item Simultaneously, compute approximations of $x$ using the given description and construct $g$ on $\bigcup_{i \leq n} A_i$ for $n=1,2,\ldots$.
\item Together, (ii) and (iii) imply that for some $n$ we will be able to verify that
$$x \in \left( \bigcup_{i \leq n} A_i \right)\degree$$
where the interior is with respect to the subset topology on $I^2$, of course.
\item When such an $n$ has been identified, we already know the definition of $g$ for $\bigcup_{i \leq n} A_i$ as well as the modulus of continuity of $g$ on the same set. This is now used to describe $g(x)$.
\end{itemize}
It remains to be shown that a suitable sequence of rectangles $\left( A_n \right)_{n \geq 1}$ exists. This follows from the next fact, assumed without proof for now\footnote{Later we shall deduce the fact from the existence of a Kleene tree.}.
\begin{lemma}[see \cite{Miller2004}, for example]
There exist computable sequences of rational numbers $(a_n)$ and $(b_n)$ in the interval $I=[0,1]$ such that the intervals $J_n=[a_n,b_n]$ have the following properties.
\begin{enumerate}[(i)]
\item If $n \neq m$ then $\left| J_n \cap J_m \right] \leq 1$.
\item If $a_n \neq 0$ then $a_n \in \{b_0,b_1,\ldots\}$ and if $b_n \neq 1$ then $b_n \in \{a_0,a_1,\ldots\}$.
\item $I_c \subsetneq \bigcup_n J_n$, i.e. the $J_n$ cover the computable reals in $I=[0,1]$.
\end{enumerate}
\end{lemma}
Now, let $\left( A_n \right)_{n \geq 1}$ be any computable enumeration of the $J_k \times J_\ell$. This completes the proof of the lemma, and the example.

\section{Baigger's example for the Polish school}
\label{sec:Bai}

\begin{figure}
\centering
\psset{unit=3cm}
\begin{pspicture}(-0.25,-0.25)(1.25,1.25)
\psaxes(0,0)(1,1)
\psframe[dimen=inner,linewidth=0pt](0,0)(1,1)
\psframe[fillcolor=white,fillstyle=solid,linestyle=none](0.1,0.1)(0.5,0.5)
\psframe[fillcolor=gray,fillstyle=solid,linestyle=none](0,0)(0.25,0.25)
\psframe[fillcolor=gray,fillstyle=solid,linestyle=none](0.75,0.75)(0.85,0.85)
\psframe[fillcolor=gray,fillstyle=solid,linestyle=none](0,0.75)(0.25,0.85)
\psframe[fillcolor=gray,fillstyle=solid,linestyle=none](0.75,0)(0.85,0.25)
\psframe[fillcolor=white,fillstyle=solid,linestyle=none](1.05,0.25)(1.25,0.75)
\dotnode(0.25,0.25){x}
\uput[u](0.275,0.225){$t_2$}
\uput[r](1,0.5){$C_2$}
\end{pspicture}
\quad
\begin{pspicture}(-0.25,-0.25)(1.25,1.25)
\psaxes(0,0)(1,1)
\psframe[dimen=inner,linewidth=0pt](0,0)(1,1)
\psframe[fillcolor=white,fillstyle=solid,linestyle=none](0.1,0.1)(0.5,0.5)
\psframe[fillcolor=gray,fillstyle=solid,linestyle=solid,linecolor=gray](0,0)(0.25,0.25)
\psframe[fillcolor=gray,fillstyle=solid,linestyle=solid,linecolor=gray](0.75,0.75)(0.85,0.85)
\psframe[fillcolor=gray,fillstyle=solid,linestyle=solid,linecolor=gray](0,0.75)(0.25,0.85)
\psframe[fillcolor=gray,fillstyle=solid,linestyle=solid,linecolor=gray](0.75,0)(0.85,0.25)
\psframe[fillcolor=gray,fillstyle=solid,linestyle=solid,linecolor=gray](0.55,0.75)(0.6,0.85)
\psframe[fillcolor=gray,fillstyle=solid,linestyle=solid,linecolor=gray](0.55,0)(0.6,0.25)
\psframe[fillcolor=gray,fillstyle=solid,linestyle=solid,linecolor=gray](0.55,0.55)(0.6,0.6)
\psframe[fillcolor=gray,fillstyle=solid,linestyle=solid,linecolor=gray](0.75,0.55)(0.85,0.6)
\psframe[fillcolor=gray,fillstyle=solid,linestyle=solid,linecolor=gray](0,0.55)(0.25,0.6)
\psframe[fillcolor=gray,fillstyle=solid,linestyle=solid,linecolor=gray](0.3,0.75)(0.4,0.85)
\psframe[fillcolor=gray,fillstyle=solid,linestyle=solid,linecolor=gray](0.3,0)(0.4,0.25)
\psframe[fillcolor=gray,fillstyle=solid,linestyle=solid,linecolor=gray](0.3,0.55)(0.4,0.6)
\psframe[fillcolor=gray,fillstyle=solid,linestyle=solid,linecolor=gray](0.3,0.3)(0.4,0.4)
\psframe[fillcolor=gray,fillstyle=solid,linestyle=solid,linecolor=gray](0.75,0.3)(0.85,0.4)
\psframe[fillcolor=gray,fillstyle=solid,linestyle=solid,linecolor=gray](0,0.3)(0.25,0.4)
\psframe[fillcolor=gray,fillstyle=solid,linestyle=solid,linecolor=gray](0.55,0.3)(0.6,0.4)
\psframe[fillcolor=gray,fillstyle=solid,linestyle=solid,linecolor=gray](0.25,0.75)(0.3,0.85)
\psframe[fillcolor=gray,fillstyle=solid,linestyle=solid,linecolor=gray](0.25,0)(0.3,0.25)
\psframe[fillcolor=gray,fillstyle=solid,linestyle=solid,linecolor=gray](0.25,0.55)(0.3,0.6)
\psframe[fillcolor=gray,fillstyle=solid,linestyle=solid,linecolor=gray](0.25,0.3)(0.3,0.4)
\psframe[fillcolor=gray,fillstyle=solid,linestyle=solid,linecolor=gray](0.25,0.25)(0.3,0.3)
\psframe[fillcolor=gray,fillstyle=solid,linestyle=solid,linecolor=gray](0.75,0.25)(0.85,0.3)
\psframe[fillcolor=gray,fillstyle=solid,linestyle=solid,linecolor=gray](0,0.25)(0.25,0.3)
\psframe[fillcolor=gray,fillstyle=solid,linestyle=solid,linecolor=gray](0.3,0.25)(0.4,0.3)
\psframe[fillcolor=gray,fillstyle=solid,linestyle=solid,linecolor=gray](0.55,0.25)(0.6,0.3)
\psframe[fillcolor=white,fillstyle=solid,linestyle=none](1.05,0.25)(1.25,0.75)
\dotnode(0.4,0.4){x}
\uput[u](0.425,0.375){$t_5$}
\uput[r](1,0.5){$C_5$}
\end{pspicture}
\caption{The ``target points'' $t_n$}
\label{fig:Bai}
\end{figure}
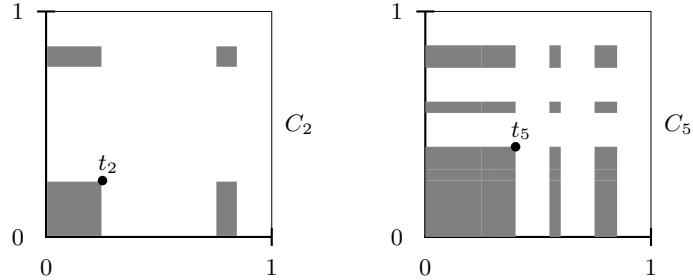

Let $a$ be any non-computable point in $I^2$. Consider the function $f$ which moves each point half-way to $a$,
$$f(x) = x + \frac{1}{2}\left( a - x \right)$$
and has a single fixed point, namely $a$ itself. The function $f$ is continuous and defined on all of $I^2$ and has no computable fixed point. Nevertheless, this is not really interesting since 
\begin{itemize}
\item the fixed point $a$ has no reasonable description---since it is itself not computable; and therefore
\item the function $f$ has no reasonable description---it is not computable in any sense.
\end{itemize}
One would like to see a function which is computable, defined (and therefore continuous) on all of $I^2$ and yet avoids fixing any of the computable points $I_c^2$. The following example, having appeared in \cite{Baigger} and in \cite{WR1999a}, modifies the construction of Orevkov to produce a computable $f$ defined on all of $I^2$ having no computable fixed point. One uses the intervals $J_n=[a_n,b_n]$ of Orevkov's example and sets 
$$C_n = \bigcup_{k,\ell \leq n} J_k \times J_\ell$$
after which one defines $f$ progressively, using the sets $C_n$. The points
$$t_n = \left(v_n, v_n \right) $$
where
$$v_n = \min_{x \in I} \left\{ x ~| ~ (x,x) \not\in C_n \right\}$$
are used as ``target point'' at each stage of the construction, as in Figure \ref{fig:Bai}. Note that
$$v = \lim_{n\rightarrow\infty} v_n$$
is not a computable number and $(v,v)$ will be one of the fixed points of $f$.

\begin{definition}
For any $W \subseteq I^2$ we define
$$W^{\blacksquare\varepsilon} = \left\{ x \in W ~\left|~ d\left(x, \delta W \setminus \delta I^2\right) \geq \varepsilon \right. \right\}$$
and
$$W^{\square\varepsilon} = \left\{ x \in W ~\left|~ d\left(x, \delta W \setminus \delta I^2\right) = \varepsilon \right. \right\}.$$
\end{definition}
\vfill
One can define $f_n$ such that 
\begin{enumerate}
\item $f_n$ moves every point in the \textit{interior} of $C_n^{\blacksquare 2^{-n}}$ but is the identity outside the set, and  is computable;
\item $f_{n+1}$ agrees with $f_n$ on $C_n^{\blacksquare 2^{-n}\cdot\frac{3}{2}}$ and therefore
\item $f = \lim_{n \rightarrow \infty} f_n$ is computable.
\end{enumerate}
Every computable point eventually lies in some
$$C_n^{\blacksquare 2^{-n}\cdot\frac{3}{2}} \quad \subset \quad \left( C_n^{\blacksquare 2^{-n}} \right)\degree$$
and is therefore moved by $f$. Clearly $f(I^2) \subseteq I^2$ and $f$ will be as required. In fact, $f$ has no fixed point in 
$$ \bigcup_n C_n = \bigcup_{k,\ell \geq 1} J_k \times J_\ell.$$
Also, $f$ has no isolated fixed point---its fixed points all occur on horizontal and vertical lines spanning the height and breadth of the unit square. Further details of the construction appear in Appendix A. The construction cannot be applied in the one-dimensional case because it is impossible to effect a change of direction by continuous rotation.

\section{BFPT and the König lemma}

In reverse mathematics it is known that in RCA$_0$, the system of recursive comprehension and $\Sigma_1^0$ -induction, the weak König lemma, WKL$_0$, is equivalent to the Brouwer FPT \cite{ShiojiTanaka}.
\begin{lemma}[WKL$_0$, Kőnig]
Every infinite binary tree has an infinite branch.
\end{lemma}
The König lemma does not have a direct computable counterpart.
\begin{theorem}[Kleene \cite{kleene_recursive_1952}]
There exists an infinite binary tree, all the computable paths of which are finite.
\end{theorem}
The relation of the Kleene tree to the Baigger counterexample is reviewed in this section. The discussion is informal and attempts only to give the essential ideas that have been revealed by the approach of reverse mathematics. In RCA$_0$, the weak König lemma WKL$_0$ has been shown to be equivalent to a number of other results in elementary analysis, such as the fact that any continuous function on a compact interval is also uniformly continuous \cite{simpson_which_1984}. WKL$_0$ and RCA$_0$ can, furthermore, be used to prove Gödel's incompleteness theorem for a countable language \cite{simpson_subsystems_1999}.

\subsection{From Baigger $f$ to Kleene tree}

Let $f$ be a computable function, as in the Baigger example, mapping $I^2$ to itself---with no computable fixed point. The following auxiliary result will be used to construct the Kleene tree.
\begin{lemma}
Let a computable $g:I^2 \rightarrow [0,1]$ be given. Then there exists a Turing-computable $h: \NZ^9 \rightarrow \NZ^2$ such that for any $(n_1,n_2,\ldots,n_8,k)$ with
$$0 \leq \frac{n_1}{n_2} \leq \frac{n_3}{n_4} \leq 1 \quad \mbox{and} \quad 0 \leq \frac{n_5}{n_6} \leq \frac{n_7}{n_8} \leq 1$$
we have $h: (n_1,n_2,\ldots,n_8,k) \mapsto (m_1,m_2)$ with $m_1 \leq m_2$ where
$$\frac{m_1}{m_2} \leq \min g \left( \left[ \frac{n_1}{n_2}, \frac{n_3}{n_4} \right] \times \left[ \frac{n_5}{n_6}, \frac{n_7}{n_8} \right]  \right) \leq \frac{m_1}{m_2} + \frac{1}{k}.$$
\end{lemma}
Let $g=||f(x)-x||$ and let $h$ be as in the lemma. Note that $g(x)=0$ if and only if $x$ is a fixed point of $f$. We shall use only the essential consequences that
\begin{itemize}
\item $g(x)>0$ for all computable $x$; and
\item there exists a (non-computable) $x_0$ such that $g(x_0)=0$.
\end{itemize}
As usual, $\{0,1\}^*$ denotes the set of finite binary sequences and $ab$ is the concatenation of $a$ and $b$.
\begin{definition}
A \textit{binary tree} is a function $t:\{0,1\}^* \rightarrow \{0,1\}$ such that 
$$t(ab)=0 \quad \mbox{for all $b$} \quad \mbox{whenever} \quad t(a)=0.$$
An infinite branch of a tree $t$ is an infinite binary sequence, on all of which finite initial segments $t$ takes the values $1$.
\end{definition}
The tree is \textit{computable} whenever the function $t$ is Turing-computable and a \textit{computable branch} is a computable binary sequence which is an infinite branch. Define the Kleene tree as follows. Let
$$t(i_1\ldots i_n) = \prod_{m=1}^{n} s(i_1\ldots i_m)$$
where $s$ is a function taking values in $\{0,1\}$. This definition of $t$ ensures that $t$ is in fact a tree and if $s$ is computable, $t$ will be a computable tree. The function $s$ will use $h$ to estimate whether $g$ gets close to zero on a specific square and if $g$ has been bounded away from zero on the square, that branch of the tree will terminate.

Define $s: \{0,1\}^* \rightarrow \{0,1\}$ for all sequences $i_1j_1\ldots i_nj_n$ of even length by
$$s(i_1j_1\ldots i_nj_n)=\chi_{\{0\}}\left(\frac{m_1}{m_2}\right)$$
where
$$(m_1,m_2) = h\left(i_1\ldots i_n,2^n,i_1\ldots i_n+1,2^n,j_1\ldots j_n,2^n,j_1\ldots j_n+1,2^n,n \right)$$
and binary strings have been interpreted as the natural numbers which they represent. Let $s$ take the value $1$ on sequences of odd length.

The tree $t$ defined in this way is obviously computable. It remains to show that $t$ is
\begin{itemize}
\item infinite; and
\item has no infinite computable branch.
\end{itemize}
Let $x_0$ be any point where $g(x_0)=0$. Then there exist infinite sequences $(i_n)$ and $(j_n)$ such that
$$x_0 \in \left[ \frac{i_1\ldots i_n}{2^n}, \frac{j_1\ldots j_n}{2^n} \right] \times \left[ \frac{i_1\ldots i_n+1}{2^n}, \frac{j_1\ldots j_n+1}{2^n} \right]\quad \mbox{for all $n$}$$
and therefore, for all $n$, $s(i_1j_1\ldots i_nj_n)=1$ and so $t(i_1j_1\ldots i_nj_n)=1$ which proves the existence of an infinite branch, hence that the tree  $t$ is infinite.

\par
Suppose that $t$ had an infinite computable branch. The branch would correspond to a decreasing chain of closed squares, the intersection of which would be non-empty. Let $x_1$ be a point in the intersection. Since, by construction of the tree, $g(x_1) \leq \frac{1}{n}$ for all $n$, $g(x_1)=0$ and hence $x_1$ would be a fixed point of $f$. However, by the construction---the branch being computable---the point $x_1$ would also be computable, contradicting the fact that $f$ has not computable fixed point. Therefore the tree $t$ has no infinite computable branch.

\subsection{From Kleene tree to Baigger $f$}

%\frametitle{From Kleene tree to Baigger $f$}
Suppose we are given a computable tree $t$ with no infinite computable branch. This tree can be used to construct a sequence of closed intervals with a computable sequence of end-points, covering all the computable real numbers in the unit interval and for which the corresponding open intervals are pair-wise disjoint.

Using the computable function $t$, one can enumerate all of the maximal finite branches of the tree. Say,
$$b(n) =b_1(n) \ldots b_{\lambda(n)}(n)$$
and set
\begin{eqnarray*}
J_{n,1} & = & \left[ \frac{b_1(n) \ldots b_{\lambda(n)}}{2^{\lambda(n)}}, \frac{b_1(n) \ldots b_{\lambda(n)}+\frac{1}{2}}{2^{\lambda(n)}} \right] \\
J_{n,m} & = & \left[ \frac{b_1(n) \ldots b_{\lambda(n)}+2^{-m+1}}{2^{\lambda(n)}}, \frac{b_1(n) \ldots b_{\lambda(n)}+2^{-m}}{2^{\lambda(n)}} \right] \quad \mbox{for $m\geq 2$.}
\end{eqnarray*}
It remains to show that the union of the intervals $J_{n,m}$ covers \textit{all} the computable points $I_c$ but not all of the unit interval $I$. It is easy to see that
\begin{itemize}
\item for every computable $x\in I_c$ there exists a computable binary sequence $(x_n)$ such that
$$\frac{x_1\ldots x_n}{2^n} \leq x < \frac{x_1\ldots x_n+1}{2^n} \quad \mbox{for all $n$}$$
and since $t$ has no infinite computable branch $t(x_1\ldots x_\ell)=0$ for some least $\ell$, in which case $x\in \cup_m J_{n,m}$ where $b(n)=x_1\ldots x_\ell$;
\item if $(x_n)$ is an infinite branch of $t$ then, since it is not computable, for all $w$ we have \mbox{$x_1x_2\ldots \neq w1111\ldots$} and therefore
$$\lim_n \frac{x_1\ldots x_n+1}{2^n} \not\in \bigcup_m J_{\ell,m}$$
for every $\ell$.
\end{itemize}
The Baigger example $f$ can now be constructed using the intervals $J_{n,m}$ and by that construction one obtains a computable $f$ with no computable fixed point, as required.

\section{Conclusion}

The existence of the Kleene tree can quite easily be derived from the impossibility of ensuring the existence of a computable fixed point for a computable function (in both Russian and Polish senses), in two dimensions (or higher). The ingenuous constructions of Orevkov and Baigger provide a way of defining a computable function with no computable fixed point from the set of intervals derived from the Kleene tree, in a constructive manner. This correspondence is, perhaps, more attractive for the ``working mathematician'' than the elegant derivation of the result in reverse mathematics. In one dimension, any computable $f: I \rightarrow I$ does have a computable point $x \in I_c$ such that $f(x)=x$, which can be seen by fairly straight-forward reduction \textit{ad absurdum} from the assumption that this is not the case.

\newpage
\bibliographystyle{splncs}
%\bibliography{BIB_v01_php_BrouwerCon}

%\newpage

\section*{Appendix A: details of the construction in Section \ref{sec:Bai}}

The constructions should guarantee that at each stage, the function $f_n$ moves every point of
$$D_n = \left( C_n^{\blacksquare 2^{-n}} \setminus C_n^{\blacksquare 2^{-n}\cdot\frac{5}{4}} \right)\degree$$
in the direction of $t_n$ by an amount proportional to its distance to $C_n^{\square 2^{-n}}$.
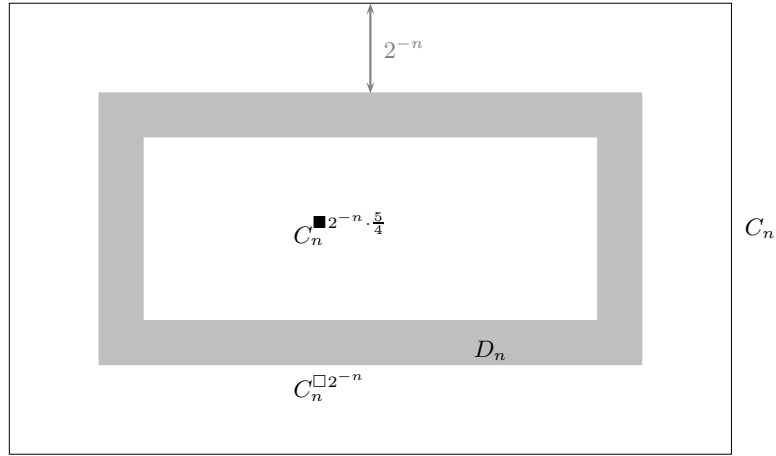
\begin{figure}
\centering
\psset{unit=6cm}
\begin{pspicture}[shift=-0.5](-0.25,0)(2,1)
\psframe[dimen=inner,linewidth=0pt](0,0)(1.6,1)
\pspolygon[linecolor=lightgray,fillstyle=solid,fillcolor=lightgray](0.2,0.2)(1.4,0.2)(1.4,0.8)(0.2,0.8)
\pspolygon[linecolor=white,fillstyle=solid,fillcolor=white](0.3,0.3)(1.3,0.3)(1.3,0.7)(0.3,0.7)
\uput[r](1.6,0.5){$C_n$}
\uput[r](0.6,0.15){$C_n^{\square 2^{-n}}$}
\uput[r](0.6,0.5){$C_n^{\blacksquare 2^{-n}\cdot\frac{5}{4}}$}
\uput[r](1.0,0.23){$D_n$}
\psline[linecolor=gray]{<->}(0.8,0.8)(0.8,1)
\uput[r](0.8,0.9){\color{gray}$2^{-n}$}
\end{pspicture}
\caption{Sets used in the construction}
\label{fig:Ba2}
\end{figure}
The construction of $f_1$ with this property is trivial. We proceed to construct $f_{n+1}$ from $f_n$.
\begin{enumerate}[(i)]
\item Extend and modify $f_n$ to $C_{n+1}^{\blacksquare 2^{-n}}$ so that every point $x$ of 
$$\left( C_{n+1}^{\blacksquare 2^{-n}} \setminus C_{n+1}^{\blacksquare 2^{-n}\cdot\frac{5}{4}} \right)\degree$$
is moved in the direction of $t_n$ by an amount proportional to $d\left(x,C_{n+1}^{\square 2^{-n}}\right)$.
%This is a simply patching of regions on which $f_n$ has already been defined.
\item Modify the resulting function so that each point in
$$C_{n+1}^{\blacksquare 2^{-n}} \setminus C_{n+1}^{\blacksquare 2^{-n}\cdot\frac{9}{8}}$$
is mapped a non-negative amount proportional to its distance to $C_{n+1}^{\square 2^{-(n+1)}}$ in the direction of $t_n$.
\item By rotation of the direction of the mapping, extend the function to $C_{n+1}^{\blacksquare 2^{-(n+1)}}$ such that every point $x$ of
$$D_{n+1} = \left( C_{n+1}^{\blacksquare 2^{-(n+1)}} \setminus C_{n+1}^{\blacksquare 2^{-(n+1)}\cdot\frac{5}{4}} \right)\degree$$
is mapped in the direction of $t_{n+1}$ by an amount proportional to $d\left(x,C_{n+1}^{\square 2^{-(n+1)}}\right)$.
\end{enumerate}
The final step is the only one in which we use the fact that we are working in two dimensions as this step requires the continuous (computable) rotation of a vector in the direction of $t_n$ to a vector in the direction of $t_{n+1}$. 

A construction is given explicitly in \cite{Baigger} but it should be clear from the preceding that it can be done in many different ways. The important part of the proof is that the construction is, at each stage, extended at the boundary to ``look right'' from the outside. This ensures that, eventually every point is in fact moved towards one of a sequence of points that converge to the non-computable fixed point $(v,v)$ on the diagonal. The Baigger construction is a somewhat delicate construction of a function that is in fact computable but that---somehow---mimics a simple mapping of every point in $I^2$ in the direction of $(v,v)$.

\end{document}